\documentclass[oneside,english]{amsart}
\usepackage[T1]{fontenc}
\usepackage[latin9]{inputenc}
\usepackage{babel}
\usepackage{amstext}
\usepackage{amsthm}
\usepackage[unicode=true,pdfusetitle,
 bookmarks=true,bookmarksnumbered=false,bookmarksopen=false,
 breaklinks=false,pdfborder={0 0 1},backref=false,colorlinks=false]
 {hyperref}

\makeatletter
\numberwithin{equation}{section}
 \theoremstyle{definition}
 \newtheorem*{defn*}{\protect\definitionname}
  \theoremstyle{plain}
  \newtheorem{thm}{\protect\theoremname}[section]
  \theoremstyle{plain}
  \newtheorem{cor}{\protect\corollaryname}[section]
  \theoremstyle{definition}
  \newtheorem{example}{\protect\examplename}[section]
  \theoremstyle{remark}
  \newtheorem{rem}{\protect\remarkname}[section]

\makeatother

  \providecommand{\definitionname}{Definition}
  \providecommand{\examplename}{Example}
  \providecommand{\remarkname}{Remark}
\providecommand{\corollaryname}{Corollary}
\providecommand{\theoremname}{Theorem}

\begin{document}

\title[$q$-analogues of series expansions for certain constants]{Two $q$-summation formulas and $q$-analogues of series expansions
for certain constants }

\author{Bing He}

\address{School of Mathematics and Statistics, Central South University\\
 Changsha 410083, Hunan, People's Republic of China}

\email{yuhe001@foxmail.com; yuhelingyun@foxmail.com}

\author{Hongcun Zhai}

\address{Department of Mathematics, Luoyang Normal University\\
 Luoyang 471934, People's Republic of China}

\email{zhai$_{-}$hc@163.com}
\begin{abstract}
From two $q$-summation formulas we deduce certain series expansion
formulas involving the $q$-gamma function. With these formulas we
can give $q$-analogues of series expansions for certain constants.

\end{abstract}

\thanks{The first author is the corresponding author. }


\maketitle

\section{Introduction}

Throughout this paper we always assume that $|q|<1.$ The $q$-gamma
function $\Gamma_{q}(x)$, first introduced by Thomae and later by
Jackson, is defined as \cite[p. 20]{GR}
\begin{equation}
\Gamma_{q}(x)=\dfrac{(q;q)_{\infty}}{(q^{x};q)_{\infty}}(1-q)^{1-x},\label{eq:1-3}
\end{equation}
where $(z;q)_{\infty}$ is given by
\[
(z;q)_{\infty}=\prod_{n=0}^{\infty}(1-zq^{n}).
\]
When $q\rightarrow1,$ the $q$-gamma function reduces to the classical
gamma function $\Gamma(x)$ which is defined by \cite{AAR}: for $\mathrm{Re}\;x>0,$
\[
\Gamma(x)=\int_{0}^{\infty}t^{x-1}e^{-t}dt.
\]
From the definition of the $q$-gamma function we know that
\begin{equation}
\dfrac{(q^{x};q)_{n}}{(1-q)^{n}}=\dfrac{\Gamma_{q}(x+n)}{\Gamma_{q}(x)},\label{eq:1-1}
\end{equation}
where $(z;q)_{n}$ is the $q$-shifted factorial given by 
\[
(z;q)_{0}=1,\;(z;q)_{n}=\prod_{k=0}^{n-1}(1-zq^{k})\;\mathrm{for}\;n\geq1.
\]
We now extend the definition of $(q^{x};q)_{n}$ to any complex $\alpha.$ 
\begin{defn*}
For any complex $\alpha,$ we define the general $q$-shifted factorial
by
\begin{equation}
(q^{x};q)_{\alpha}=\dfrac{\Gamma_{q}(x+\alpha)}{\Gamma_{q}(x)}(1-q)^{\alpha}\label{eq:1-2}
\end{equation}

For brevity, we denote $\dfrac{(q^{x};q)_{\alpha}}{(1-q)^{\alpha}}$
by $(x|q)_{\alpha},$ namely, $(x|q)_{\alpha}=\dfrac{\Gamma_{q}(x+\alpha)}{\Gamma_{q}(x)}.$
For any non-negative integer $n,$ we have
\[
(x|q)_{n}=\prod_{k=0}^{n-1}[x+k]_{q}
\]
and
\[
(x|q)_{-n}=\dfrac{\Gamma_{q}(x-n)}{\Gamma_{q}(x)}=\dfrac{(1-q)^{n}}{(q^{x-n};q)_{n}},
\]
 where $[z]_{q}$ is the $q$-integer defined by 
\[
[z]_{q}=\dfrac{1-q^{z}}{1-q}.
\]
In particular, 
\[
(x|q)_{0}=1,\;(x|q)_{1}=[x]_{q},\;(x|q)_{-1}=\frac{1}{[x-1]_{q}}.
\]
\end{defn*}
Gosper in \cite{G} introduced $q$-analogues of $\sin x$ and $\pi:$
\[
\sin_{q}(\pi x):=q^{(x-1\text{/2})^{2}}\frac{(q^{2-2x};q^{2})_{\infty}(q^{2x};q^{2})_{\infty}}{(q;q^{2})_{\infty}^{2}}
\]
and 
\[
\pi_{q}:=(1-q^{2})q^{1/4}\dfrac{(q^{2};q^{2})_{\infty}^{2}}{(q;q^{2})_{\infty}^{2}}.
\]
They satisfy the following relations:
\[
\lim_{q\rightarrow1}\sin_{q}x=\sin x,\quad\lim_{q\rightarrow1}\pi_{q}=\pi
\]
and 
\begin{equation}
\Gamma_{q^{2}}(x)\Gamma_{q^{2}}(1-x)=\dfrac{\pi_{q}}{\sin_{q}(\pi x)}q^{x(x-1)}.\label{eq:1-4}
\end{equation}
When $q\rightarrow1,$ the last identity reduces to the Euler reflection
formula \cite[(1.2.1)]{AAR}:
\[
\Gamma(x)\Gamma(1-x)=\dfrac{\pi}{\sin(\pi x)}.
\]

Ramanujan \cite{R} recorded without proof 17 series expansions for
$1/\pi$, among which, the proof of the first three was briefly sketched
in \cite{R1}. The first complete proof of all 17 formulas was found
by the Borwein brothers \cite{BB}. D.V. Chudnovsky and G.V. Chudnovsky
\cite{CC} proved some of the Ramanujan\textquoteright s series representations
for $1/\pi$ independently and established new series as well. The
readers can refer to the paper \cite{BBC} for the history of the
Ramanujan-type series for $1/\pi.$ Recently, using certain properties
of the general rising shifted factorial and the gamma function, Liu
in \cite{L1,L2} supplied many series expansion formula for $1/\pi.$
$q$-Analogues of two Ramanujan-type series for $1/\pi$ were established
by Guo and Liu \cite{GL} using $q$-WZ pairs and some basic hypergeometric
identities. 

Our motivation for the present work emanates from \cite{GL,L1,L2}.
In this paper we shall deduce from two $q$-summation formulas certain
series expansion formulas involving the $q$-gamma function. These
formulas allow us to give $q$-analogues of series expansions for
certain constants. These series expansion formulas are as follows.
\begin{thm}
\label{t1} For any complex number $\alpha$ and $\mathrm{Re}(c-a-b)>0$
we have
\[
\begin{gathered}\sum_{n=0}^{\infty}\dfrac{(\alpha|q^{2})_{a+n}(1-\alpha|q^{2})_{b+n}}{[n]_{q^{2}}!\Gamma_{q^{2}}(c+n+1)}q^{2(c-a-b)}{}^{n}\\
=\dfrac{(\alpha|q^{2})_{a}(1-\alpha|q^{2})_{b}\Gamma_{q^{2}}(c-a-b)}{(1-\alpha|q^{2})_{c-a}(\alpha|q^{2})_{c-b}}q^{-\alpha(\alpha-1)}\cdot\dfrac{\sin_{q}(\pi\alpha)}{\pi_{q}},
\end{gathered}
\]
where $[n]_{q}!$ is given by 
\[
[0]_{q}!=1,\quad[n]_{q}!=\prod_{k=1}^{n}[k]_{q}\;\mathrm{for}\;n\geq1.
\]
\end{thm}
\begin{thm}
\label{t} For $\mathrm{Re}(a+b+c+d+1+\alpha-\beta-\gamma-\delta)>0$
we have
\[
\begin{gathered}\sum_{n=0}^{\infty}\frac{(1-q^{4n+2a+2\alpha})(\alpha|q^{2})_{a+n}(\beta|q^{2})_{n-b}(\gamma|q^{2})_{n-c}(\delta|q^{2})_{n-d}}{(1-q^{2})[n]_{q^{2}}!(1+\alpha-\beta|q^{2})_{a+b+n}(1+\alpha-\gamma|q^{2})_{a+c+n}(1+\alpha-\delta|q^{2})_{a+d+n}}q^{An}\\
=\frac{\Gamma_{q^{2}}(1+\alpha-\beta)\Gamma_{q^{2}}(1+\alpha-\gamma)\Gamma_{q^{2}}(1+\alpha-\delta)\Gamma_{q^{2}}(2+\alpha-\beta-\gamma-\delta)}{\Gamma_{q^{2}}(\alpha)\Gamma_{q^{2}}(1+\alpha-\beta-\gamma)\Gamma_{q^{2}}(1+\alpha-\beta-\delta)\Gamma_{q^{2}}(1+\alpha-\gamma-\delta)}\\
\times\frac{(\beta|q^{2})_{-b}(\gamma|q^{2})_{-c}(\delta|q^{2})_{-d}(2+\alpha-\beta-\gamma-\delta|q^{2})_{a+b+c+d-1}}{(1+\alpha-\beta-\gamma|q^{2})_{a+b+c}(1+\alpha-\beta-\delta|q^{2})_{a+b+d}(1+\alpha-\gamma-\delta|q^{2})_{a+c+d}},
\end{gathered}
\]
where $A=2(a+b+c+d+1+\alpha-\beta-\gamma-\delta).$
\end{thm}
The next section is devoted to our proof of Theorems \ref{t1} and
\ref{t}. In Section \ref{sec:R} we deduce $q$-analogues of certain
series expansions for $1/\pi.$ In the last section several $q$-analogues
of series expansions for $\pi^{2}$ are also obtained. 

\section{Proof of Theorems \ref{t1} and \ref{t}}

\noindent{\it Proof of Theorem \ref{t1}.} Recall from \cite[(1.5.1)]{GR}
the $q$-Gauss summation formula:
\[
\sum_{n=0}^{\infty}\dfrac{(a;q)_{n}(b;q)_{n}}{(q;q)_{n}(c;q)_{n}}\big(c/ab\big)^{n}=\dfrac{(c/a;q)_{\infty}(c/b;q)_{\infty}}{(c;q)_{\infty}(c/ab;q)_{\infty}},\quad|c/ab|<1.
\]
Making the substitutions: $q\rightarrow q^{2},\;a\rightarrow q^{2a},\;b\rightarrow q^{2b},\;c\rightarrow q^{2c}$
in the above identity and using (\ref{eq:1-3}) and (\ref{eq:1-1})
we have
\[
\sum_{n=0}^{\infty}\dfrac{\Gamma_{q^{2}}(a+n)\Gamma_{q^{2}}(b+n)}{[n]_{q^{2}}!\Gamma_{q^{2}}(c+n)}q^{2(c-a-b)}{}^{n}=\dfrac{\Gamma_{q^{2}}(a)\Gamma_{q^{2}}(b)\Gamma_{q^{2}}(c-a-b)}{\Gamma_{q^{2}}(c-a)\Gamma_{q^{2}}(c-b)}.
\]
Replacing $a,\;b,\;c$ by $a+\alpha,\;b+1-\alpha,\;c+1$ respectively
in the above formula we get
\begin{equation}
\begin{gathered}\sum_{n=0}^{\infty}\dfrac{\Gamma_{q^{2}}(a+n+\alpha)\Gamma_{q^{2}}(b+n+1-\alpha)}{[n]_{q^{2}}!\Gamma_{q^{2}}(c+n+1)}q^{2(c-a-b)}{}^{n}\\
=\dfrac{\Gamma_{q^{2}}(a+\alpha)\Gamma_{q^{2}}(b+1-\alpha)\Gamma_{q^{2}}(c-a-b)}{\Gamma_{q^{2}}(c-a+1-\alpha)\Gamma_{q^{2}}(c-b+\alpha)}.
\end{gathered}
\label{eq:2-1}
\end{equation}
It follows from (\ref{eq:1-2}) that 
\[
\begin{gathered}\Gamma_{q^{2}}(a+\alpha)=(\alpha|q^{2})_{a}\Gamma_{q^{2}}(\alpha),\\
\Gamma_{q^{2}}(b+1-\alpha)=(1-\alpha|q^{2})_{b}\Gamma_{q^{2}}(1-\alpha),\\
\Gamma_{q^{2}}(a+n+\alpha)=(\alpha|q^{2})_{a+n}\Gamma_{q^{2}}(\alpha),\\
\Gamma_{q^{2}}(b+n+1-\alpha)=(1-\alpha|q^{2})_{b+n}\Gamma_{q^{2}}(1-\alpha),\\
\Gamma_{q^{2}}(c-a+1-\alpha)=(1-\alpha|q^{2})_{c-a}\Gamma_{q^{2}}(1-\alpha),\\
\Gamma_{q^{2}}(c-b+\alpha)=(\alpha|q^{2})_{c-b}\Gamma_{q^{2}}(\alpha).
\end{gathered}
\]
Substituting these formulas into (\ref{eq:2-1}) and simplifying we
arrive at
\[
\sum_{n=0}^{\infty}\dfrac{(\alpha|q^{2})_{a+n}(1-\alpha|q^{2})_{b+n}}{[n]_{q^{2}}!\Gamma_{q^{2}}(c+n+1)}q^{2(c-a-b)}{}^{n}=\dfrac{(\alpha|q^{2})_{a}(1-\alpha|q^{2})_{b}\Gamma_{q^{2}}(c-a-b)}{(1-\alpha|q^{2})_{c-a}(\alpha|q^{2})_{c-b}\Gamma_{q^{2}}(\alpha)\Gamma_{q^{2}}(1-\alpha)}.
\]
From this identity and (\ref{eq:1-4}) we can deduce the result readily.
This completes the proof of Theorem \ref{t1}. \qed

\noindent{\it Proof of Theorem \ref{t}.} Recall the following summation
formula for the basic hypergeometric series \cite[(2.7.1)]{GR}:
\begin{equation}
_{6}\phi_{5}\left(\begin{matrix}a,qa^{\frac{1}{2}},-qa^{\frac{1}{2}},b,c,d\\
a^{\frac{1}{2}},-a^{\frac{1}{2}},aq/b,aq/c,aq/d
\end{matrix};q,\dfrac{aq}{bcd}\right)=\dfrac{(aq,aq/bc,aq/bd,aq/cd;q)_{\infty}}{(aq/b,aq/c,aq/d,aq/bcd;q)_{\infty}},\label{eq:1}
\end{equation}
where $\bigg|\dfrac{aq}{bcd}\bigg|<1$ and $_{6}\phi_{5}$ is the
basic hypergeometric series given by 
\[
_{6}\phi_{5}\left(\begin{matrix}a_{1},a_{2},a_{3},a_{4},a_{5},a_{6}\\
b_{1},b_{2},b_{3},b_{4},b_{5}
\end{matrix};q,z\right)=\sum_{n=0}^{\infty}\frac{(a_{1},a_{2},a_{3},a_{4},a_{5},a_{6};q)_{n}}{(q,b_{1},b_{2},b_{3},b_{4},b_{5};q)_{n}}z^{n}.
\]
Replacing $(q,a,b,c,d)$ by $(q^{2},q^{2a},q^{2b},q^{2c},q^{2d})$
in (\ref{eq:1}) and employing (\ref{eq:1-1}) and (\ref{eq:1-2})
we have
\begin{equation}
\begin{gathered}\sum_{n=0}^{\infty}\frac{(1-q^{4n+2a})\Gamma_{q^{2}}(a+n)\Gamma_{q^{2}}(b+n)\Gamma_{q^{2}}(c+n)\Gamma_{q^{2}}(d+n)}{(1-q^{2})[n]_{q^{2}}!\Gamma_{q^{2}}(1+a-b+n)\Gamma_{q^{2}}(1+a-c+n)\Gamma_{q^{2}}(1+a-d+n)}q^{2n(1+a-b-c-d)}\\
=\frac{\Gamma_{q^{2}}(b)\Gamma_{q^{2}}(c)\Gamma_{q^{2}}(d)\Gamma_{q^{2}}(1+a-b-c-d)}{\Gamma_{q^{2}}(1+a-b-c)\Gamma_{q^{2}}(1+a-b-d)\Gamma_{q^{2}}(1+a-c-d)}.
\end{gathered}
\label{eq:2}
\end{equation}
It follows from (\ref{eq:1-3}) that 

\[
\begin{gathered}\Gamma_{q^{2}}(a+n+\alpha)=(\alpha|q^{2})_{a+n}\Gamma_{q^{2}}(\alpha),\Gamma_{q^{2}}(n-b+\beta)=(\beta|q^{2})_{n-b}\Gamma_{q^{2}}(\beta),\\
\Gamma_{q^{2}}(n-c+\gamma)=(\gamma|q^{2})_{n-c}\Gamma_{q^{2}}(\gamma),\Gamma_{q^{2}}(n-d+\delta)=(\delta|q^{2})_{n-d}\Gamma_{q^{2}}(\delta),\\
\Gamma_{q^{2}}(\beta-b)=(\beta|q^{2})_{-b}\Gamma_{q^{2}}(\beta),\Gamma_{q^{2}}(\gamma-c)=(\gamma|q^{2})_{-c}\Gamma_{q^{2}}(\gamma),\Gamma_{q^{2}}(\delta-d)=(\delta|q^{2})_{-d}\Gamma_{q^{2}}(\delta),\\
\Gamma_{q^{2}}(a+b+n+1+\alpha-\beta)=(1+\alpha-\beta|q^{2})_{a+b+n}\Gamma_{q^{2}}(1+\alpha-\beta),\\
\Gamma_{q^{2}}(a+c+n+1+\alpha-\gamma)=(1+\alpha-\gamma|q^{2})_{a+c+n}\Gamma_{q^{2}}(1+\alpha-\gamma),\\
\Gamma_{q^{2}}(a+d+n+1+\alpha-\delta)=(1+\alpha-\delta|q^{2})_{a+d+n}\Gamma_{q^{2}}(1+\alpha-\delta),\\
\Gamma_{q^{2}}(a+b+c+1+\alpha-\beta-\gamma)=(1+\alpha-\beta-\gamma|q^{2})_{a+b+c}\Gamma_{q^{2}}(1+\alpha-\beta-\gamma),\\
\Gamma_{q^{2}}(a+b+d+1+\alpha-\beta-\delta)=(1+\alpha-\beta-\delta|q^{2})_{a+b+d}\Gamma_{q^{2}}(1+\alpha-\beta-\delta),\\
\Gamma_{q^{2}}(a+c+d+1+\alpha-\gamma-\delta)=(1+\alpha-\gamma-\delta|q^{2})_{a+c+d}\Gamma_{q^{2}}(1+\alpha-\gamma-\delta)
\end{gathered}
\]
and
\[
\begin{gathered}\Gamma_{q^{2}}(a+b+c+d+1+\alpha-\beta-\gamma-\delta)\\
=(2+\alpha-\beta-\gamma-\delta|q^{2})_{a+b+c+d-1}\Gamma_{q^{2}}(2+\alpha-\beta-\gamma-\delta).
\end{gathered}
\]
Making the substitutions: $a\rightarrow a+\alpha,\:b\rightarrow\beta-b,\:c\rightarrow\gamma-c,\:d\rightarrow\delta-d$
in (\ref{eq:2}) and then substituting the above identities into the
resulting equation we can easily deduce the result. This finishes
the proof of Theorem \ref{t}. \qed

\section{\label{sec:R} $q$-Analogues of certain series expansions for $1/\pi$}

In this section we employ Theorems \ref{t1} and \ref{t} to deduce
$q$-analogues of certain series expansions for $1/\pi.$

Setting $\alpha=\dfrac{1}{2}$ in Theorem \ref{t1} and using the
fact $\sin_{q}\dfrac{\pi}{2}=1$ we get 
\begin{thm}
\label{t3} For $\mathrm{Re}(c-a-b)>0$ we have
\[
\sum_{n=0}^{\infty}\dfrac{(1/2|q^{2})_{a+n}(1/2|q^{2})_{b+n}}{[n]_{q^{2}}!\Gamma_{q^{2}}(c+n+1)}q^{2(c-a-b)}{}^{n}=\dfrac{(1/2|q^{2})_{a}(1/2|q^{2})_{b}\Gamma_{q^{2}}(c-a-b)}{(1/2|q^{2})_{c-a}(1/2|q^{2})_{c-b}}\cdot\dfrac{q^{1/4}}{\pi_{q}}.
\]
\end{thm}
We put $a=b=0$ and $c=l$ in Theorem \ref{t3} to arrive at
\begin{cor}
If $l$ is positive integer, then
\[
\sum_{n=0}^{\infty}\dfrac{(1/2|q^{2})_{n}^{2}}{[n]_{q^{2}}!(l|q^{2})_{n+1}}q{}^{2ln}=\dfrac{q^{1/4}}{\pi_{q}(1/2|q^{2})_{l}^{2}}.
\]
\end{cor}
\begin{example}
($l=1$) We have

\[
\sum_{n=0}^{\infty}\dfrac{(1/2|q^{2})_{n}^{2}}{[n]_{q^{2}}![n+1]_{q^{2}}!}q{}^{2n}=\dfrac{(1+q)^{2}q^{1/4}}{\pi_{q}}.
\]
\end{example}
This expansion for $1/\pi_{q}$ is a $q$-analogue of the series for
$1/\pi$ \cite[p. 174]{Gl}:
\[
\sum_{n=0}^{\infty}\frac{(1/2)_{n}^{2}}{n!(n+1)!}=\frac{4}{\pi},
\]
where $(1/2)_{n}$ is the shifted factorial given by 
\[
(1/2)_{0}=1,\:(1/2)_{n}=\prod_{k=0}^{n-1}(1/2+k)\:\mathrm{for}\:n\geq1.
\]
Actually, this expansion for $1/\pi_{q}$ was also obtained by Guo
\cite[(1.8)]{Guo}.
\begin{example}
($l=2$) We have
\[
\sum_{n=0}^{\infty}\dfrac{(1/2|q^{2})_{n}^{2}}{[n]_{q^{2}}![n+2]_{q^{2}}!}q{}^{4n}=\dfrac{(1+q)^{4}q^{1/4}}{\pi_{q}(1+q+q^{2})^{2}}.
\]
\end{example}
This series expansion for $1/\pi_{q}$ can be considered as a $q$-analogue
of the series for $1/\pi:$
\[
\sum_{n=0}^{\infty}\frac{(1/2)_{n}^{2}}{n!(n+2)!}=\frac{16}{9\pi}.
\]

We set $a=b=-1$ and $c=l$ in Theorem \ref{t3} to deduce 
\begin{cor}
If $l$ is a non-negative integer, then
\[
\begin{gathered}q^{2}(1+q)^{2}[l+1]_{q^{2}}+q^{2l+4}+[l+1]_{q^{2}}!\sum_{n=1}^{\infty}\dfrac{(1/2|q^{2})_{n}^{2}}{[n+1]_{q^{2}}![l+n+1]_{q^{2}}!}q^{2(l+2)(}{}^{n+1)}\\
=\dfrac{(1+q)^{2}[l+1]_{q^{2}}!^{2}}{\pi_{q}(1/2|q^{2})_{l+1}^{2}}q^{9/4}.
\end{gathered}
\]
\end{cor}
\begin{example}
($l=0$) We have 
\[
q^{2}(1+q)^{2}+q^{4}+\sum_{n=1}^{\infty}\dfrac{(1/2|q^{2})_{n}^{2}}{[n+1]_{q^{2}}!^{2}}q^{4}{}^{n+4}=\dfrac{(1+q)^{4}}{\pi_{q}}q^{9/4}.
\]
\end{example}
This series expansion for $1/\pi_{q}$ can be regarded as a $q$-analogue
of the series for $1/\pi$ \cite[p. 174]{Gl}:
\[
5+\sum_{n=1}^{\infty}\frac{(1/2)_{n}^{2}}{(n+1)!^{2}}=\frac{16}{\pi}.
\]
\begin{example}
($l=1$) We have 
\[
q^{2}(1+q)^{2}(1+q^{2})+q^{6}+(1+q^{2})\sum_{n=1}^{\infty}\dfrac{(1/2|q^{2})_{n}^{2}}{[n+1]_{q^{2}}![n+2]_{q^{2}}!}q^{6}{}^{n+6}=\dfrac{(1+q)^{6}(1+q^{2})^{2}}{\pi_{q}(1+q+q^{2})^{2}}q^{9/4}.
\]
\end{example}
This expansion for $1/\pi_{q}$ is also a $q$-analogue of the series
for $1/\pi:$
\[
9+2\sum_{n=1}^{\infty}\frac{(1/2)_{n}^{2}}{(n+1)!(n+2)!}=\frac{256}{9\pi}.
\]
\begin{rem}
Besides those formulas displayed in Theorem \ref{t3} and its consequences,
we can give some other new series expansions for $1/\pi_{q}$ with
the change of $\alpha.$ We shall not display them out one by one
in this paper.
\end{rem}
\begin{thm}
\label{t3-1} For $\mathrm{Re}(a+b+c+d)>0$ we have
\[
\begin{gathered}\sum_{n=0}^{\infty}\frac{(1-q^{4n+2a+1})(1/2|q^{2})_{a+n}(1/2|q^{2})_{n-b}(1/3|q^{2})_{n-c}(2/3|q^{2})_{n-d}}{(1-q^{2})[n]_{q^{2}}!(1|q^{2})_{a+b+n}(7/6|q^{2})_{a+c+n}(5/6|q^{2})_{a+d+n}}q^{2(a+b+c+d)n}\\
=\frac{(1/2|q^{2})_{-b}(1/3|q^{2})_{-c}(2/3|q^{2})_{-d}(1|q^{2})_{a+b+c+d-1}}{(1/3|q^{2})_{a+b+d}(2/3|q^{2})_{a+b+c}(1/2|q^{2})_{a+c+d}}\cdot\frac{[1/6]_{q^{2}}(q^{4/3},q^{2/3};q^{2})_{\infty}q^{1/4}}{(q^{1/3},q^{5/3};q^{2})_{\infty}\pi_{q}}.
\end{gathered}
\]
\end{thm}
\noindent{\it Proof.} It follows from (\ref{eq:1-4}) that 
\begin{align}
\Gamma_{q^{2}}^{2}(1/2) & =\pi_{q}q^{-1/4},\label{eq:3-1}\\
\Gamma_{q^{2}}(1/3)\Gamma_{q^{2}}(2/3) & =\dfrac{\pi_{q}}{\sin_{q}(\pi/3)}q^{-2/9},\nonumber \\
\Gamma_{q^{2}}(7/6)\Gamma_{q^{2}}(5/6) & =[1/6]_{q^{2}}\Gamma_{q^{2}}(1/6)\Gamma_{q^{2}}(5/6)\nonumber \\
 & =\dfrac{\pi_{q}}{\sin_{q}(\pi/6)}[1/6]_{q^{2}}q^{-5/36}.\nonumber 
\end{align}
Then, by the definition of $\sin_{q},$
\begin{equation}
\frac{\Gamma_{q^{2}}(7/6)\Gamma_{q^{2}}(5/6)}{\Gamma_{q^{2}}(1/3)\Gamma_{q^{2}}(2/3)}=\frac{\sin_{q}(\pi/3)}{\sin_{q}(\pi/6)}[1/6]_{q^{2}}q^{1/12}=\frac{(q^{4/3},q^{2/3};q^{2})_{\infty}[1/6]_{q^{2}}}{(q^{1/3},q^{5/3};q^{2})_{\infty}}.\label{eq:3-2}
\end{equation}
Therefore, the result follows easily by setting $(\alpha,\beta,\gamma,\delta)=(1/2,1/2,1/3,2/3)$
in Theorem \ref{t} and applying the identities $\Gamma_{q}(1)=1,$
(\ref{eq:3-1}) and (\ref{eq:3-2}). \qed

Taking $(a,b,c,d)=(1,0,0,0)$ in Theorem \ref{t3-1} we can get
\begin{example}
We have
\[
\begin{gathered}\sum_{n=0}^{\infty}\frac{(1-q^{4n+3})(1-q^{2n+1})(1/2|q^{2})_{n}^{2}(1/3|q^{2})_{n}(2/3|q^{2})_{n}}{(1-q^{2})(1-q^{2n+2})([n]_{q^{2}}!)^{2}(7/6|q^{2})_{1+n}(5/6|q^{2})_{1+n}}q^{2n}\\
=\frac{[1/6]_{q^{2}}(q^{4/3},q^{2/3};q^{2})_{\infty}q^{1/4}}{[1/3]_{q^{2}}[2/3]_{q^{2}}[1/2]_{q^{2}}(q^{1/3},q^{5/3};q^{2})_{\infty}\pi_{q}}.
\end{gathered}
\]
\end{example}
This series expansion for $1/\pi_{q}$ can be regarded as a $q$-analogue
of the series for $1/\pi:$
\[
\sum_{n=0}^{\infty}\frac{(4n+3)(2n+1)(1/2)_{n}^{2}(1/3)_{n}(2/3)_{n}}{(n+1)(6n+1)(6n+5)(6n+7)(n!)^{2}(1/6)_{n}(5/6)_{n}}=\frac{\sqrt{3}}{6\pi}.
\]

Putting $(a,b,c,d)=(0,0,0,1)$ in Theorem \ref{t3-1} we can deduce
that
\begin{example}
We have
\[
\begin{gathered}\frac{q^{2/3}}{(1+q)[1/3]_{q^{2}}[5/6]_{q^{2}}}-\sum_{n=1}^{\infty}\frac{(1-q^{4n+1})(1/2|q^{2})_{n}^{2}(1/3|q^{2})_{n}(2/3|q^{2})_{n-1}}{(1-q^{2})([n]_{q^{2}}!)^{2}(7/6|q^{2})_{n}(5/6|q^{2})_{1+n}}q^{2n}\\
=\frac{[1/6]_{q^{2}}}{[1/3]_{q^{2}}^{2}[1/2]_{q^{2}}}\cdot\frac{(q^{4/3},q^{2/3};q^{2})_{\infty}q^{11/12}}{(q^{1/3},q^{5/3};q^{2})_{\infty}\pi_{q}}.
\end{gathered}
\]
\end{example}
This series expansion for $1/\pi_{q}$ can be considered as a $q$-analogue
of the series for $1/\pi:$
\[
1-\frac{5}{18}\sum_{n=1}^{\infty}\frac{(4n+1)(1/2)_{n}^{2}(1/3)_{n}(2/3)_{n-1}}{(n!)^{2}(7/6)_{n}(5/6)_{1+n}}=\frac{5}{\sqrt{3}\pi}.
\]

\section{$q$-Analogues of series expansions for $\pi^{2}$}

In this section we use Theorem \ref{t} to give $q$-analogues of
some series expansions for $\pi^{2}.$
\begin{thm}
\label{t4-1} For $\mathrm{Re}(a+b+c+d-1/2)>0$ we have

\[
\begin{gathered}\sum_{n=0}^{\infty}\frac{(1-q^{4n+2a})(1|q^{2})_{a+n-1}(1/2|q^{2})_{n-b}(1/2|q^{2})_{n-c}(1/2|q^{2})_{n-d}}{(1-q^{2})[n]_{q^{2}}!(1/2|q^{2})_{a+b+n}(1/2|q^{2})_{a+c+n}(1/2|q^{2})_{a+d+n}}q^{2(a+b+c+d)n-n}\\
=\frac{\pi_{q}^{2}(1/2|q^{2})_{-b}(1/2|q^{2})_{-c}(1/2|q^{2})_{-d}(1/2|q^{2})_{a+b+c+d-1}}{(1|q^{2})_{a+b+c-1}(1|q^{2})_{a+b+d-1}(1|q^{2})_{a+c+d-1}q^{1/2}}
\end{gathered}
\]
\end{thm}
\noindent{\it Proof.} It can be dedeuced from $\Gamma_{q}(x+1)=[x]_{q}\Gamma_{q}(x)$
and Theorem \ref{t} that 

\[
\begin{gathered}\sum_{n=0}^{\infty}\frac{(1-q^{4n+2a+2\alpha})(\alpha+1|q^{2})_{a+n-1}(\beta|q^{2})_{n-b}(\gamma|q^{2})_{n-c}(\delta|q^{2})_{n-d}}{(1-q^{2})[n]_{q^{2}}!(1+\alpha-\beta|q^{2})_{a+b+n}(1+\alpha-\gamma|q^{2})_{a+c+n}(1+\alpha-\delta|q^{2})_{a+d+n}}q^{An}\\
=\frac{\Gamma_{q^{2}}(1+\alpha-\beta)\Gamma_{q^{2}}(1+\alpha-\gamma)\Gamma_{q^{2}}(1+\alpha-\delta)\Gamma_{q^{2}}(2+\alpha-\beta-\gamma-\delta)}{\Gamma_{q^{2}}(\alpha+1)\Gamma_{q^{2}}(2+\alpha-\beta-\gamma)\Gamma_{q^{2}}(2+\alpha-\beta-\delta)\Gamma_{q^{2}}(2+\alpha-\gamma-\delta)}\\
\times\frac{(\beta|q^{2})_{-b}(\gamma|q^{2})_{-c}(\delta|q^{2})_{-d}(2+\alpha-\beta-\gamma-\delta|q^{2})_{a+b+c+d-1}}{(2+\alpha-\beta-\gamma|q^{2})_{a+b+c-1}(2+\alpha-\beta-\delta|q^{2})_{a+b+d-1}(2+\alpha-\gamma-\delta|q^{2})_{a+c+d-1}}.
\end{gathered}
\]
Then the result follows readily from by setting $(\alpha,\beta,\gamma,\delta)=(0,1/2,1/2,1/2)$
in the above identity and applying the identities $\Gamma_{q}(1)=1$
and (\ref{eq:3-1}). \qed

Taking $(a,b,c,d)=(1,0,0,0)$ in Theorem \ref{t4-1} we can obtain
\begin{example}
We have
\[
\sum_{n=0}^{\infty}\frac{(1+q^{2n+1})q^{n}}{(1-q^{2n+1})^{2}}=\frac{\pi_{q}^{2}}{(1-q^{2})^{2}q^{1/2}}.
\]
\end{example}
This series expansion for $\pi_{q}^{2}$ can be regarded as a $q$-analogue
of the series for $\pi^{2}:$
\[
\sum_{n=0}^{\infty}\frac{1}{(2n+1)^{2}}=\frac{\pi^{2}}{8}.
\]
Actually, this expansion for $\pi_{q}^{2}$ has been obtained by Sun
\cite[(1.2)]{S}.

Setting $(a,b,c,d)=(1,1,1,0)$ in Theorem \ref{t4-1} we can derive
\begin{example}
We have
\[
\sum_{n=0}^{\infty}\frac{(1+q^{2n+1})q^{5n}}{(1-q^{2n-1})^{2}(1-q^{2n+1})^{2}(1-q^{2n+3})^{2}}=\frac{\pi_{q}^{2}(1+q+q^{2})q^{3/2}}{(1+q^{2})(1-q^{2})^{6}}.
\]
\end{example}
This series expansion for $\pi_{q}^{2}$ can also be considered as
a $q$-analogue of the series for $\pi^{2}:$
\[
\sum_{n=0}^{\infty}\frac{1}{(2n-1)^{2}(2n+1)^{2}(2n+3)^{2}}=\frac{3\pi^{2}}{256}.
\]

Putting $(a,b,c,d)=(1,1,1,1)$ in Theorem \ref{t4-1} we can deduce 
\begin{example}
We have
\[
\begin{gathered}\frac{(1+q)q^{3}}{(1-q)^{5}(1-q^{3})^{3}}-\sum_{n=1}^{\infty}\frac{(1+q^{2n+1})q^{7n}}{(1-q^{2n-1})^{3}(1-q^{2n+1})^{2}(1-q^{2n+3})^{3}}\\
=\frac{\pi_{q}^{2}(1+q+q^{2})(1+q+q^{2}+q^{3}+q^{4})q^{5/2}}{(1+q^{2})^{3}(1-q^{2})^{8}}.
\end{gathered}
\]
\end{example}
This series expansion for $\pi_{q}^{2}$ is also a $q$-analogue of
the series for $\pi^{2}:$
\[
\frac{1}{27}-\sum_{n=1}^{\infty}\frac{1}{(2n-1)^{3}(2n+1)^{2}(2n+3)^{3}}=\frac{15\pi^{2}}{4096}.
\]
\begin{rem}
Besides those formulas displayed in Theorems \ref{t3-1} and \ref{t4-1}
and their consequences, we can give a general series expansion for
$1/\pi_{q}^{2}$ by taking $(\alpha,\beta,\gamma,\delta)=(1/2,1/2,1/2,1/2)$
in Theorem \ref{t}, from which many series expansions for $1/\pi_{q}^{2}$
can be deduced. We shall not display them out one by one in this work. 
\end{rem}

\section*{Acknowledgements}

 The first author was partially supported by the National Natural
Science Foundation of China (Grant No. 11801451). The second author
was supported by the National Natural Science Foundation of China
(Grant No. 11371184) and the Natural Science Foundation of Henan Province
(Grant No. 162300410086, 2016B259, 172102410069).


\begin{thebibliography}{10}
\bibitem{AAR}G.E. Andrews, R. Askey and R. Roy, Special Functions,
Encyclopedia of Mathematics and its Applications, Vol. 71. Cambridge
University Press, Cambridge, 1999.

\bibitem{BBC}N.D. Baruah, B.C. Berndt and H.H. Chan, Ramanujan\textquoteright s
series for $1/\pi:$ A survey, Amer. Math. Monthly 116 (2009), 567\textendash 587.

\bibitem{BB}J.M. Borwein and P.B. Borwein, Pi and the AGM. Wiley,
New York, 1987.

\bibitem{CC}D.V. Chudnovsky, G.V. Chudnovsky, Approximation and complex
multiplication according to Ramanujan, in: G.E. Andrews, R.A. Askey,
B.C. Berndt, K.G. Ramanathan, R.A. Rankin (Eds.), Ramanujan Revisited,
Academic Press, Boston, 1988, pp. 375\textendash 472.

\bibitem{GR}G. Gasper and M. Rahman, Basic Hypergeometric Series,
Cambridge University Press, Cambridge, 1990.

\bibitem{Gl}J.W.L. Glaisher, On series for $1/\pi$ and $1/\pi^{2},$
Quart. J. Pure Appl. Math. 37 (1905), 173\textendash 198.

\bibitem{G}R.W. Gosper, Experiments and discoveries in $q$-trigonometry,
in: F.G. Garvan, M.E.H. Ismail (Eds.), Symbolic Computation, Number
Theory, Special Functions, Physics and Combinatorics, Kluwer, Dordrecht,
Netherlands, 2001, pp.79\textendash 105.

\bibitem{Guo}V.J.W. Guo, A $q$-Analogue of the (I.2) Supercongruence
of Van Hamme, International Journal of Number Theory, doi: 10.1142/S1793042118501701.

\bibitem{GL}V.J.W. Guo and J.-C. Liu, $q$-analogues of two Ramanujan-type
formulas for $1/\pi,$ arXiv: 1802.01944v2.

\bibitem{L1}Z.-G. Liu, A summation formula and Ramanujan type series,
J. Math. Anal. Appl. 389(2) (2012), 1059\textendash 1065.

\bibitem{L2}Z.-G. Liu, Gauss summation and Ramanujan-type series
for $1/\pi.$ Int. J. Number Theory, 8(2) (2012), 289\textendash 297.

\bibitem{R1}S. Ramanujan, Collected Papers, Cambridge University
Press, Cambridge, 1927, reprinted by Chelsea, New York, 1962, reprinted
by the American Mathematical Society, Providence, RI, 2000.

\bibitem{R}S. Ramanujan, Modular equations and approximations to
$\pi,$ Quart. J. Pure Appl. Math. 45 (1914), 350\textendash 372.

\bibitem{S}Z.-W. Sun, Two $q$-analogues of Euler's formula $\zeta(2)=\pi^{2}/6,$
arXiv: 1802.01473v3.
\end{thebibliography}
\end{document}